\documentclass[journal]{IEEEtran}

\usepackage{ifpdf}

\usepackage[ruled]{algorithm2e}
\usepackage{color}
\usepackage{url}
\ifpdf
\usepackage[pdftex]{graphicx}
\else
\usepackage[dvips]{graphicx}
\fi
\usepackage[caption=false,font=footnotesize]{subfig}
\usepackage{comment}
\usepackage{amsmath,amsfonts} 
\usepackage[noadjust]{cite}
\usepackage{nomencl}
\usepackage{bm}
\newtheorem{theorem}{\textbf{Theorem}}

\newtheorem{lemma}{\textbf{Lemma}}

\usepackage[latin1]{inputenc}
\usepackage{tikz}
\usetikzlibrary{shapes,arrows}
\usepackage{xcolor}

\definecolor{ggray}{gray}{.9}

\newcommand{\ud}[1]{_\mathrm{#1}}
\newcommand{\up}[1]{^\mathrm{#1}}
\newcommand\Tstrut{\rule{0pt}{2.6ex}}         
\newcommand\Bstrut{\rule[-0.9ex]{0pt}{0pt}}   

\tikzstyle{decision} = [diamond, draw, fill=blue!20, 
    text width=4.5em, text badly centered, node distance=3cm, inner sep=0pt]
\tikzstyle{block} = [rectangle, draw, fill=blue!20, 
    text width=5em, text centered, rounded corners, minimum height=4em]
\tikzstyle{line} = [draw, -latex']
\tikzstyle{cloud} = [draw, ellipse,fill=red!20, node distance=5cm,
    minimum height=2em]

\IEEEoverridecommandlockouts
\begin{document}

\title{Optimal Battery Participation \\in Frequency Regulation Markets}

\author{
    Bolun Xu,~\IEEEmembership{Student Member,~IEEE,}
    Yuanyuan Shi,~\IEEEmembership{Student Member,~IEEE,}\\
    Daniel S. Kirschen~\IEEEmembership{Fellow,~IEEE,}
    and Baosen Zhang~\IEEEmembership{Member,~IEEE}
\thanks{The authors are with the Department of Electrical Engineering,
University of Washington,
Seattle, Washington 98125,
(e-mail:\{xubolun, yyshi, kirschen, zhangbao\}@uw.edu). This work has been supported in part by the University of Washington Clean Energy Institute. }
}

\maketitle
\makenomenclature

\begin{abstract}
Battery participants in performance-based frequency regulation markets must consider the cost of battery aging in their operating strategies to maximize market profits. In this paper we solve this problem by proposing an optimal control policy and an optimal bidding policy based on realistic market settings and an accurate battery aging model. The proposed control policy has a threshold structure and achieves near-optimal performance with respect to an offline controller that has complete future information. The proposed bidding policy considers the optimal control policy to maximize market profits while satisfying the market performance requirement through a chance-constraint. It factors the value of performance and supports a  trade-off between higher profits and a lower risk of violating performance requirements. We demonstrate the optimality of both policies using simulations. A case study based on the PJM regulation market shows that our approach is effective at maximizing operating profits.
\end{abstract}


\IEEEpeerreviewmaketitle

\section*{Nomenclature}
\addcontentsline{toc}{section}{Nomenclature}

\subsection{Parameters and Variables}

\begin{IEEEdescription}[\IEEEusemathlabelsep\IEEEsetlabelwidth{$V_1,V_2,V_3$}]
\item[$B$] Battery energy storage power rating in MW
\item[$b_t$] Battery dispatch power during $t$ in MW
\item[$\bm{b}$] The set of all battery dispatch power $\bm{b}=\{b_t\}$
\item[$C$] Regulation capacity in MW
\item[$\overline{C}$] Maximum regulation capacity 
\item[$E$] Battery energy storage capacity in MWh
\item[$\overline{E}, \underline{E}$] Upper and lower energy limit for the battery energy storage in MWh
\item[$\overline{E}^g_t, \underline{E}^g_t$] Upper and lower energy limit enforced by the control policy $g$ during interval $t$ in MWh
\item[$e_t$] Battery energy level during $t$ in MWh
\item[$\bm{e}$] The set of all battery energy levels $\bm{e}=\{e_t\}$
\item[$e\up{max}_t,e\up{min}_t$] Maximum and minimum energy level up to interval $t$ in MWh
\item[$g$] Control policy for the battery energy storage
\item[$i$] Index of identified cycles
\item[$j$] Index of regulation capacity bid segments
\item[$J$] Number of total bid segments
\item[$M$] Duration of one dispatch interval in hours
\item[$R$] Battery cell replacement price in \$/MWh
\item[$r_t$] Normalized regulation instruction during $t$
\item[$\bm{r}$] The set of all normalized regulation instructions $\bm{r}=\{r_t\}$
\item[$t$] Dispatch time interval index
\item[$T$] Number of intervals during a market period
\item[$u_i$] The depth of cycle $i$
\item[$\bm{u}$] Set of all cycle depths
\item[]
\item[$\gamma$] Relative regulation capacity in hours
\item[$\delta$] Coefficient in the performance index model
\item[$\lambda$] Regulation market clearing price in \$/MW
\item[$\eta$] Battery charge and discharge efficiency
\item[$\pi$] Performance penalty price in \$/MWh
\item[$\mu_{\lambda}$] Expected market clearing price in \$/MW
\item[$\mu_{r}$] Expected total regulation instruction
\item[$\rho\up{min}$] Minimum performance index requirement
\item[$\xi$] Performance chance constraint confidence level

\subsection{Functions}
\item[$A(\cdot)$] Cycle aging cost function
\item[$J(\cdot)$] Regulation market profit expectation
\item[$P(\cdot)$] Performance index calculation function
\item[$\Pi(\cdot)$] Regulation market profit
\item[$\Phi(\cdot)$] Cycle depth aging stress function
\item[$\varphi(\cdot)$] Derivative of $\Phi(\cdot)$
\end{IEEEdescription}

\section{Introduction}

The share of battery energy storage (BES) in the frequency regulation markets is increasing rapidly~\cite{irena}. In the PJM market, the BES capacity has increased from zero in 2005 to over 280~MW in 2017, making up 41\% of its regulation procurement capacity~\cite{pjm2016}.
This growth rate was made possible by the rapid decrease in battery cell manufacturing cost~\cite{nykvist2015rapidly}, the increase in renewable penetration~\cite{qiu2016stochastic}, as well as changes in market rules that have lowered the barrier for batteries to provide frequency regulation. Following the Federal Energy Regulatory Commission (FERC) Orders~755~\cite{FERC755} and ~784~\cite{FERC784}, all Independent System Operators (ISO) and Regional Transmission Organizations (RTO) in the U.S. have implemented pay-for-performance regulation markets, and improved the automatic generation control (AGC) framework to account for the state of charge (SoC) constraint of energy storage in their regulation dispatch~\cite{xu2016comparison}. 
In these revised markets, BES participants are rewarded for their fast and accurate regulation response. They are also able to offer a very small energy capacity (15 minutes energy-to-power ratio at minimum) to provide regulation products specially designed for energy storage. 

These regulation market reforms have encouraged numerous BES projects targeted solely at the provision of regulation. These BES units often behave as price-takers by submitting offers at zero price into the market and accepting any clearing prices. If the market prices for regulation are high, these participants can earn considerable profits with naive control strategies. However, since the amount of regulation capacity needed is limited, a market can easily become saturated by too many price-taking participants. Operational evidence shows that the PJM RegD market became saturated in 2016 and its market clearing prices have dropped by two-thirds since 2014~\cite{xu2014bess,pjm2016}, making BES operation hardly profitable with naive bidding and operating strategies. 


BES participants must employ more advanced operating and bidding strategies to secure operational profits against dropping market prices and the cost of battery degradation. Previous studies have assessed how providing regulation with batteries affects their aging characteristic under simple SoC control strategies~\cite{mercier2009optimizing, borsche2013power,watson2017comparing,xu2014bess,tu2017maximizing,byrne2012estimating,zhang2017data, swierczynski2014selection, stroe2015degradation, koller2015review, white2011using, stroe2016degradation}. However, few studies have actively incorporated battery aging as part of regulation operation or bidding optimization objectives. References \cite{zhang2016profit, cui2015optimal} take into account the battery lifespan in regulation control optimization, using an aging model that is too simple to reflect properly the complex battery cycle aging mechanisms. The results in \cite{he2015optimal} incorporate the battery aging cost into regulation bidding strategies, but the proposed method does not optimize real-time operations. 

Since regulation instructions are highly stochastic and occur at a very high time resolution, accounting for them explicitly in an optimization problem is computationally challenging. It is thus crucial to reduce the formulation complexity using statistical and analytical derivations, and to jointly optimize the bidding and real-time control strategies. In this work, we incorporate battery cycle aging into a real-time control policy that ensures profitability under any market clearing result. In addition, we introduce a bidding policy based on the proposed control policy that maximizes profits while satisfying the market performance requirement. The main contributions of this paper are as follows:
\begin{itemize}
\item We formulate battery bidding and operation in a pay-for-performance frequency regulation market as a chance-constrained profit maximization problem that includes an accurate battery cycle aging model. 
\item We apply a real-time control policy to optimize the battery regulation response by balancing the battery cycle aging cost and the regulation mismatch penalty. This policy achieves near-optimal control results and applies to any regulation market designs.
\item We prove that under the proposed control policy, the optimal battery regulation problem is reduced to finding the highest capacity possible while ensuring the required regulation performance. We therefore propose a bidding policy that ensures regulation performance while maximizing the operational profit by factoring the value of performance in the frequency regulation market.
\item Numerical examples and case studies using real market data are given to demonstrate the  effectiveness of the proposed policies.
\end{itemize}
The rest of this paper is organized as follows. Section~II describes the formulation of the optimization problem. Section~III introduces the proposed control and bidding policies. Section~IV describes and discusses the simulation results. Section~V concludes the paper.

\section{Problem Formulation}

\subsection{Frequency Regulation Market Model}

We consider a pay-for-performance regulation market in which a participant is rewarded based on the regulation capacity it provides, as well as on how much it is instructed to alter dispatch set-points by the system operator, and how accurately it follows the regulation instructions. Most system operators calculate a performance index as the relative error between a participant's regulation response and the regulation instructions~\cite{xu2016comparison}. PJM uses a more complicate calculation method~\cite{pjm_manual12}.

Pay-for-performance regulation markets have two-part offer and payment designs. A participant submits a regulation capacity offer and a regulation mileage offer, and then the system operator unifies the two offer prices into a single modified offer based on the participant's historical performance index and the expected dispatch regulation mileage. 
In the ex-post market settlement, the capacity payment is calculated using the capacity clearing price and the assigned regulation capacity. The mileage payment is calculated using the mileage clearing price and the instructed mileage. Depending on the market design, the performance index is used to penalize the entire regulation payment, or only the mileage payment. All pay-for-performance regulation markets also have a minimum performance requirement, in which a participant must reach a certain performance index to be eligible for receiving regulation payment, and must maintain a satisfactory performance history to be qualified to participate in the regulation market~\cite{xu2014bess}.

The objective of a regulation market participant is to maximize its operational profit, including the expected market payment and the battery aging cost. To formulate this problem, we start by considering a market period with finite discrete dispatch intervals $t\in \{1,\dotsc, T\}$ and a single market clearing price $\lambda$ which is calculated based on the capacity clearing price, the mileage clearing prices, and the instructed mileage~\cite{xu2016comparison}. 
We assume that the participant receives regulation instructions that are proportionally scaled with respect to its regulation capacity $C$ based on a set of normalized regulation signals $\bm{r}=\{r_t\in [-1,1]\}$. The participant receives the payment $P(C\bm{r}, \bm{b})\lambda C $ from the regulation market,
where $\bm{b}=\{b_t\}$ is the participant's response to the regulation instruction, and $P(\cdot)$ is the performance index calculation function. We define $P(\cdot)$ as a linear function of the relative response error:
\begin{align}\label{eq:lin_per}
    P(C\bm{r}, \bm{b}) =1-  \frac{||C\bm{r}-\bm{b}||_1}{C||\bm{r}||_1}\delta\,.
\end{align}
where $||C\bm{r}-\bm{b}||_1$ calculates the absolute error in the response, and $C||\bm{r}||_1$ is the total amount of the instructed regulation signal. $\delta \in [0,1]$ is the part of the regulation payment that is evaluated based on the relative error of the response. This is a unified performance model that fits all market rules. The value of $\delta$ is specific to each market and can be found in each market manual. 



\subsection{BES Operation Model}
The BES operation has the following constraints
\begin{align}
    -B &\leq b_t \leq B \label{eq:bat_con1}\\
    \underline{E}&\leq e_t \leq \overline{E} \label{eq:bat_con2}\\
    e_t-e_{t-1} &= M\eta[b_t]^+ - M[-b_t]^+/\eta \label{eq:bat_con3}
\end{align}
where $b_t$ is the battery dispatch power at time $t$, $B$ is the BES power ratings, $\bm{e}=\{e_t\}$ is the battery energy level at step $t$, $M$ is the dispatch interval duration, $\eta$ is the BES single-trip charge or discharge efficiency, $\overline{E}$ and $\underline{E}$ are the BES upper and lower energy limits in MWh, and  $[x]^+=\max\{x,0\}$. \eqref{eq:bat_con1} models the BES power rating, \eqref{eq:bat_con2} models the BES storage limits, and \eqref{eq:bat_con3} models the evolution of the state of charge.

\subsection{Battery Aging Model}

The chemistry of battery aging process is typically modeled with partial differential equations~\cite{ramadesigan2012modeling, ning2004cycle}. These models have good accuracy, but are highly nonlinear and pose significant computational challenges for the optimal control problems~\cite{reniers2018improving}. Instead, many studies on grid-scale BES optimization formulate battery aging as empirical linear or quadratic functions with respect to control variables including active power, energy throughput, and SoC~\cite{ying2016stochastic,wang2017improving,megel2013maximizing,fortenbacher2017modeling,ortega2014optimal,koller2013defining,wang2016quantifying,shi2017using,nguyen2016stochastic}. These empirical models can be seamlessly incorporated into linear or quadratic programming problems and solved within tractable times. In this paper, we explore an alternative battery modeling approach by incorporating a systematic cycle counting algorithm into optimization. This approach offers an efficient way of modeling the varying aging effect from cycles of different depth.
Battery operations in frequency regulation mostly consist shallow cycles due to frequent switching between charging and discharging. These shallow cycles cause much lower aging damage per MWh of energy throughput because battery cycle aging is highly nonlinear with respect to the cycle depth. For example, a Lithium Nickel Manganese Cobalt Oxide (NMC) cell can provide ten times more lifetime energy throughput if cycled at 10\% depth (between 55\% SoC and 65\% SoC) than if cycled at 100\% depth~\cite{ecker2014calendar}. This nonlinear aging property, in which a smaller cycle depth results in a higher lifetime energy throughput, is observed in most static electrochemical batteries~\cite{ruetschi2004aging, byrne2012estimating, xu2016modeling, wang2014degradation}, and it is therefore crucial to use an accurate cycle aging model for regulation operation.

The cycle aging cost model we adopt in this study relies on the rainflow method~\cite{rychlik1987new} for cycle identification. The rainflow method decomposes a set of normalized battery SoC measurements into a combination of independent cycles and the cycle aging is calculated with respect to each identified cycle. The rainflow-based aging model therefore offers much higher accuracy in irregular battery operations, and has been extensively adopted for battery life assessment and operation optimization~\cite{he2015optimal,abdulla2016optimal,jin2018applicability}. The rainflow method decomposes a SoC series into a combination of cycles as:
\begin{align}\label{eq:age_model1}
    \bm{u} &= \mathbf{Rainflow}(\bm{e}/E)\,
\end{align}
where $\bm{u}$ is the set of all cycle depths. $\bm{e}=\{e_t\}$ is the BES energy level set, $E$ is the BES rated energy capacity. $\bm{e}/E$ is thus the normalized SoC series. After identifying cycles, a cycle depth stress function $\Phi(u)$ calculates the incremental battery cycle aging caused by a particular cycle of depth $u$. The total cycle life loss is calculated by summing the life loss from all identified cycles from the SoC series. The cycle depth stress can be obtained from lab test results~\cite{laresgoiti2015modeling,millner2010modeling,omar2014lithium} which are described in the manufacturer's warranty. Since the battery cycle aging cost prorates the battery cell replacement cost to the cycle life loss, the the battery aging cost function $A(\cdot)$ is defined as follows:
\begin{align}\label{eq:A}
    A(\bm{b}) &= ER\textstyle\sum_{i=1}^{|\bm{u}|}\Phi(u_i)
\end{align}
where $|\bm{u}|$ is the cardinality of $\bm{u}$, $R$ is the battery cell replacement cost in \$/MWh which includes the disposal cost of old cells and the purchase and installation cost for new cells. $\Phi(u)$ is the cycle depth aging stress function, and a half cycle causes half the aging stress of a full cycle. The aging cost function $A$ is written as a function of $\bm{b}$ because the battery SoC is a linear combination of the active power. For example, if the cell replacement cost of a 1~MWh NMC BES is \$300k and this cell operates 100,000 cycles at 10\% cycle depth, then a 10\% depth cycle (0.1~MWh energy throughput) costs \$3. This cell operates only 1,000 cycles at 100\% cycle depth, then a 100\% depth cycle (1~MWh energy throughput) would cost \$300.

\subsection{Optimization Problem}
If the market clearing price $\lambda$ and the regulation signal realization $\bm{r}$ are known, the participant can find the optimal regulation capacity $C$ and the optimal BES dispatch $\bm{b}$ by solving the following optimization problem
\begin{subequations}\label{eq:opt0}
\begin{align}
    \max_{C, \bm{b}}\;&
    \Pi(C,\lambda,\bm{r}, \bm{b}) := P(C\bm{r}, \bm{b})\lambda C  - A(\bm{b}) 
    \label{eq:opt0_obj}\\
    \text{s.t. }&  P(C\bm{r}, \bm{b}) \geq \rho\up{min}
    \label{eq:opt0_c1}\\
    &C\in [0,B], \text{constraints \eqref{eq:bat_con1}--\eqref{eq:bat_con3}}
    \label{eq:opt0_c2}\,
\end{align}
\end{subequations}
where the objective function \eqref{eq:opt0_obj} maximizes the operating profit calculated as the market revenue minus the aging cost, \eqref{eq:opt0_c1} enforces a minimum performance index $\rho\up{min}$ to the regulation response, \eqref{eq:opt0_c2} limits the regulation capacity $C$ to be within the battery power capacity $B$, and the battery operation must satisfy the operation constraints \eqref{eq:bat_con1}--\eqref{eq:bat_con3}.

However, this problem cannot be solved in practice because the realization of $\bm{r}$ is not known in advance. The participant must therefore decide on a policy $g$ that determines $b_t$ at time step $t$ based only on past information. The expected profit $J(\cdot)$ corresponding to a particular operational policy and regulation capacity $(g,C)$ is defined as:
\begin{align}
    J(g,C) &= \mathbb{E}\Big[\Pi(C,\lambda,\bm{r}, \bm{b}^g)\Big] \nonumber\\
    & =  \mathbb{E}\Big[\lambda C -\delta\lambda\frac{||C\bm{r}-\bm{b}^g||_1}{||\bm{r}||_1}
-A(\bm{b}^g)\Big]\,,
\end{align}
where the superscript $g$ is included to indicate the dependence on the control policy $g$. However, the exact value of $\lambda$ and $||\bm{r}||_1$ can only be known after the regulation provision, hence we approximate the formulation using the following expectations
\begin{align}
    \mu_{\lambda} &= \mathbb{E}\big[\lambda\big]\,,\quad\mu_{r} = \mathbb{E}\big[||\bm{r}||_1\big]\,.
\end{align}
To solve this problem, we transform the minimum performance requirement into a chance constraint, and rewrite  \eqref{eq:opt0} as a chance-constrained stochastic programming problem
\begin{align}\label{eq:opt1}
    \max_{g,C}\; &J(g,C) = \mu_{\lambda} C - \mathbb{E}\Big[\delta\frac{\mu_{\lambda}}{\mu_{r}}
||C\bm{r}-\bm{b}^g||_1+A(\bm{b}^g)\Big]\nonumber\\
    \text{s.t. }& \mathbf{Prob}\big[P(C\bm{r}, \bm{b}) \geq \rho\up{min}\big] \geq \xi \nonumber\\
    &g \in \mathcal{G},\; C\in [0,B] \,.
 \end{align}
where $\xi$ is the confidence level that the performance requirement will  be satisfied, and $\mathcal{G}$ is the set of all feasible operation policies $g$ that satisfy \eqref{eq:bat_con1}--\eqref{eq:bat_con3} and the causality condition.

\section{Optimal Control and Bidding Policy}\label{sec:policy}

\begin{figure}
    \centering
\begin{tikzpicture}[node distance = 3cm, auto]
    \node [cloud] (bid) {bidding policy};
    \node [cloud, right of=bid] (control) {control policy};
    \node [block, below of=bid] (market) {system operator};
    \node [block, below of=control] (dispatch) {regulation dispatch};

    \path [line,dashed] (control) -- node [above,midway] {optimal design} (bid);
    \path [line,dashed] (bid) -- node [above, rotate=90, midway] {offer curve} (market);
    \path [line,dashed] (control) -- node [above, rotate=90, midway] {response} (dispatch);
    \path [line] (market) -- node [below,midway] {regulation signal} (dispatch);
    \path [line] (market) -- node [rotate=30, near end] {clearing prices} (control);
\end{tikzpicture}
\caption{Flow chart of the regulation market and control processes. The control policy inputs the clearing prices for determining the optimal response, while the bidding policy submits offers to the market while considering the control policy.}
\label{fig:chart}
\end{figure}
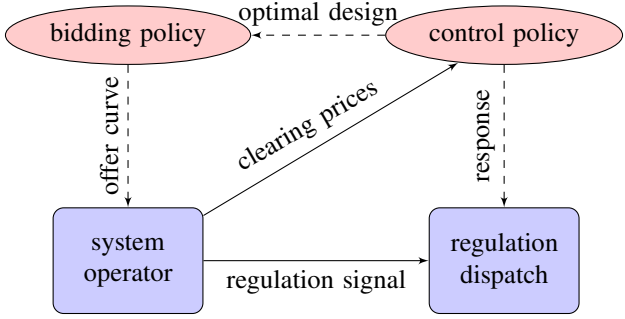

We propose an online control policy and a bidding policy that solve problem \eqref{eq:opt1}. Control and bidding in the regulation market are part of a sequential decision-making process because the regulation capacity $C$ is cleared before the actual dispatch.  We therefore work backwards by first determining the optimal online control policy $g^*$ with respect to any cleared regulation capacity (Fig.~\ref{fig:chart})
\begin{align}\label{eq:opt_g}
    g^* \in \mathrm{arg}\max_{g\in \mathcal{G}}J(g,C)          
\end{align}
We then derive the optimal bidding policy $C\up{*}$ based on the optimal control policy:
\begin{align}\label{eq:opt_C}
    C^* \in \mathrm{arg}\max_{C\in [0,B]}&\Big\{\max_{g\in \mathcal{G}}J(g,C)\Big\}\nonumber\\
    \text{s.t. }& \mathbf{Prob}\big[P(C\bm{r}, \bm{b}) \geq \rho\up{min}\big] \geq \xi \,.
\end{align}
We consider the performance chance constraint at the bidding stage because it depends primarily on the regulation capacity and the battery energy capacity (see Section~\ref{sec:cap}). 
We assume that $\bm{r}$ is energy zero-mean (including efficiency losses) because either the system operator or the participant can employ strategies for controlling the average BES SoC, and these strategies do not effect a participant's objective of minimizing its operating cost.

\subsection{Optimal Control Policy}

\begin{figure}[t]%
	\centering
	\subfloat[]{
		\includegraphics[trim = 0mm 0mm 0mm 0mm, clip, width = .95\columnwidth]{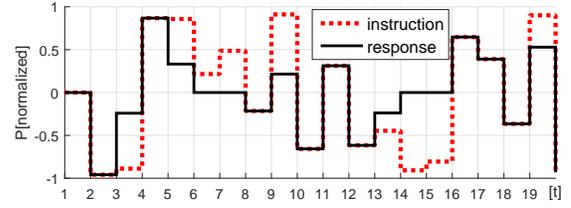}
		\label{fig:con1}%
	}
	\\
	\subfloat[Controlled vs. uncontrolled SoC profile]{
		\includegraphics[trim = 0mm 0mm 0mm 0mm, clip, width = .95\columnwidth]{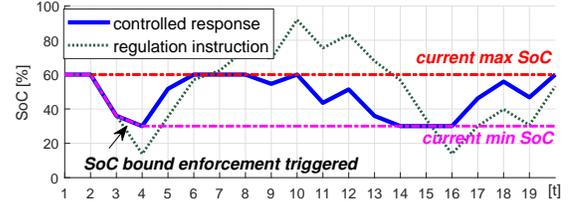}
		\label{fig:con2}%
	}
	\caption{\footnotesize Illustration of the proposed control policy. The policy keeps track of the real-time maximum and minimum SoC level. When the distance between them reaches the calculated threshold $\hat{u}$, the policy starts to constrain the response. (a) dispatch instruction vs. response; (b) Controlled vs. uncontrolled SoC profile.
	}%
	\label{fig:con}
\end{figure}

\begin{algorithm}[!htb]\label{alg:ctrl}
\SetAlgoLined
\KwResult{Determine battery dispatch points $b_t$}
 \tcp{\footnotesize calculate optimal cycle depth}
 set $\hat{u}\to\varphi^{-1}\Big(\frac{\eta^2+1}{\eta R}\pi\Big)$\;
 \tcp{\footnotesize set current max/min energy level to $e_0$}
 set $e_0 \to e\up{max}_{0}$, $e_0 \to e\up{min}_{0}$\;
 \While{$t\leq T$}{
  \tcp{\footnotesize read $e_t$ and update min/max energy level}
  set $\max\{e\up{max}_{n-1}, e_t\} \to e\up{max}_t$, $\min\{e\up{min}_{t-1}, e_t\} \to e\up{min}_t$\;
  \tcp{\footnotesize activate energy level bounds if either reached $\hat{u}$ or the state of charge constraints}

  set $\min\{\overline{E}, e\up{min}_t + \hat{u}E\} \to \overline{E}^g_t$\;
  set $\max\{\underline{E}, e\up{max}_t - \hat{u}E\} \to \underline{E}^g_t$\;
  \tcp{\footnotesize read $Cr_t$ and enforce soc bounds}
  \eIf{$Cr_t \geq 0$}{
   set $\min\Big\{\frac{1}{M\eta}(\overline{E}^g-e_t), Cr_t\Big\} \to b_t$ \;
   }{
   set $\max\Big\{\frac{\eta}{M}(\underline{E}^g-e_t), Cr_t\Big\}\to b_t$ \;
  }
  \tcp{\footnotesize wait until the next control interval}
  set $t+1\to t$\;
 }
\caption{Optimal Battery Regulation Response Policy}
\end{algorithm}

The online regulation response policy balances the cost of deviating from the regulation signal and the cycle aging cost of batteries while satisfying operating constraints. This policy takes a threshold form and achieves an optimality gap that is independent of the total number of time steps. Therefore in term of regret (consequence of decision-making under uncertainty)~\cite{bell1982regret}, this policy achieves the strongest possible result: the regret does not grow with time. The key part of the control policy is to calculate thresholds that bound the SoC of the battery as functions of the deviation penalty and degradation cost. The bound on the SoC $\hat{u}$ is given by:
\begin{align}\label{eq:pol3}
    \hat{u} = \varphi^{-1}\Big(\frac{\eta^2+1}{\eta R}\pi\Big),\quad\pi = \delta\frac{\mu_{\lambda}}{\mu_{r}M} 
\end{align}
where $\varphi^{-1}(\cdot)$ is the inverse function of the derivative of the cycle stress function $\varphi(x) = \mathrm{d} \Phi(x)/\mathrm{d} x$. Algorithm~\ref{alg:ctrl} summarizes this control policy, and Fig.~\ref{fig:con} shows an example of control based on this policy, where the battery follows the dispatch instructions until the distance between its maximum and minimum SoC reaches $\hat{u}$.

\begin{figure}
    \centering
    \includegraphics[trim = 10mm 0mm 10mm 0mm, clip, width = .9\columnwidth]{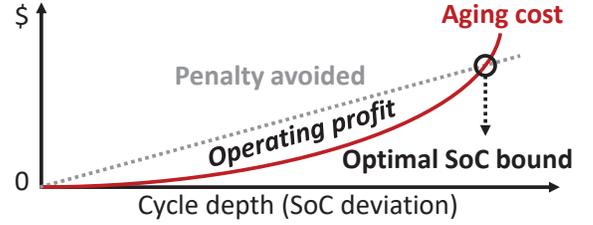}
    \caption{Illustration of the optimal cycle depth when responding to regulation instructions. During a single cycle, the participant avoids a mismatch penalty that is linear with respect to the cycle depth, and undertakes a nonlinear cycle aging cost (this relationship is explained Appendix. A). Hence there exists an optimal cycle depth that maximizes the operating profit which is the area between the two curves.}
    \label{fig:cycle_illu}
\end{figure}

Fig.~\ref{fig:cycle_illu} explains the intuition of the proposed policy. First we consider an arbitrary cycle in the regulation provision. The battery performs a cycle by following the regulation instruction and avoids a penalty cost linear to the cycle depth, but undertakes a nonlinear aging cost. Hence there exists an optimal full cycle depth $\hat{u}$ that maximizes the operating profit over a single cycle. The battery should stop following the regulation instruction once this optimal depth is reached, which is equivalent to setting a SoC bound as described in the proposed policy. In cases that the aging cost function is convex, this optimal bound can be explicitly calculated from \eqref{eq:pol3}, and the proposed policy has a bounded regret given in the following theorem.

\begin{theorem}\label{theorem1}
Suppose the battery cycle aging stress function $\Phi(\cdot)$ is strictly convex. The proposed control strategy $g(\cdot)$ has a worst-case optimality gap (regret) $\epsilon$ that is independent of the operation time duration $TN$:
\begin{align}
	&\Big[\delta\frac{\mu_{\lambda}}{\mu_{r}}||C\bm{r}-\bm{b}^*||_1-A(\bm{b}^*)\Big]-\Big[\delta\frac{\mu_{\lambda}}{\mu_{r}}||C\bm{r}-\bm{b}^g||_1-A(\bm{b}^g)\Big] \nonumber\\
	& \leq \epsilon\, \label{eq:gap}\nonumber\\
    & \text{$\forall$ $e_0 \in [\underline{E}, \overline{E}]$, $C\in [0,B]$, and}\nonumber\\
    & \text{$\forall$ sequences $\{r_t\} \in [-1,1]$, $t\in\{1,\dotsc,T\}$} 
\end{align}
where $\bm{b}^*\in \mathrm{arg}\min_{\bm{b}\in  \mathcal{B}(e_0)} \pi||C\bm{r}-\bm{b}||_1-A(\bm{b})$. 
\end{theorem}

A proof of Theorem~\ref{theorem1} can be found in \cite{shi2017optimal}, and mathematical derivations of the control policy and the optimal SoC bound $\hat{u}$ are summarized in Appendix~A. The regret exists since battery operation may contain a residue that cannot be matched into cycles. Note that the optimality of the proposed policy depends on the convexity of the cycle aging function, however the computation and implementation of the proposed policy does not require convexity and this policy applies to any batteries that have a monotonic increasing cycle stress function.



\subsection{Optimal Regulation Capacity}\label{sec:cap}

We first show that under the optimal control policy and when \emph{relaxing} the minimum performance constraint, a participant should adopt a price-taker bidding strategy, i.e., it should provide the maximum possible regulation capacity under any market prices. We then incorporate the performance chance constraint and show that the performance-constrained optimal capacity can be explicitly characterized based on statistics of historical regulation signals.

\begin{theorem}
When relaxing the minimum performance requirement, the optimal regulation capacity $C^*$ is equal to the BES power rating $B$ when using the proposed optimal response policy $g^*$.
\end{theorem}

\emph{Proof:}
See Appendix B.

The intuition for this theorem is that the proposed policy guarantees optimal operating profit under any regulation capacity and signal realization. A higher regulation capacity offers higher profit potential since market payment is capped by $\lambda C$. Therefore under optimal real-time control, a participant will never earn less profit with a higher capacity because the penalty cost never exceeds the market payment.

Stemming from Theorem~2, the solution of \eqref{eq:opt_C} is reduced to finding the maximum $C$ while satisfying the performance chance-constraint. The challenge is that the regulation signal realization is unknown when determining the regulation capacity. We must therefore characterize a probabilistic function that correlates the performance index $\rho$ with the regulation capacity $C$. In addition, since the battery regulation response $\bm{b}^g$ under the proposed policy $g^*$ is price responsive, we must also take the expected clearing price $\mu_{\lambda}$ into consideration.

We start with the performance index calculation with the regulation signal realization $\bm{r}$ and show that it can be reformed into a function with respect to $\bm{b}^g/C$, which is the normalized battery response
\begin{align}\label{eq:pi_1}
    P(C\bm{r}, \bm{b}) &=1-  \frac{||C\bm{r}-\bm{b}^g||_1}{C||\bm{r}||_1} \nonumber\\
    &= 1-  \frac{||\bm{r}-{\bm{b}^g}/{C}||_1}{||\bm{r}||_1}\,.
\end{align}
 Recall that in the proposed policy the battery energy level is constrained between $\overline{E}^g$ and $\underline{E}^g$ that are calculated from battery data and the penalty price. We substitute $\bm{b}^g/C$ into the proposed policy and show that $\bm{b}^g/C$ only depends on $\overline{E}_g/C$ and $\underline{E}_g/C$ since the battery operation constraints \eqref{eq:bat_con1}--\eqref{eq:bat_con3} are linear
\begin{align}\label{eq:pi_2}
    \frac{b^g_t}{C} = \begin{cases}
    \min\big\{{(\overline{E}^g-e_t)}/{(M\eta C)}, r_t\big\} & \text{ if $r_t\geq 0$}\\
    \max\big\{{(\underline{E}^g-e_t)}/{(M\eta C)}, r_t\big\} & \text{ if $r_t < 0$}\,.
    \end{cases}
\end{align}
We omit the effect of the battery initial energy level $e_0$ on the performance index because $\bm{r}$ has a zero mean (i.e. it is energy neutral). The performance index thus depends only on the ratio between the usable energy capacity of the battery $\overline{E}^g-\underline{E}^g$ and the regulation capacity $C$. We can therefore define this ratio as the normalized regulation energy capacity $\gamma$
\begin{align}\label{eq:pi_3}
    \gamma = \frac{\overline{E}^g - \underline{E}^g}{C}=\frac{\min\{\overline{E}-\underline{E},\hat{u}E\}}{C}
\end{align}
where we substitute $\hat{u}$ into $\overline{E}^g - \underline{E}^g$ as in Algorithm~1. Recall that $\hat{u}$ is calculated from the expected market clearing price, hence we represent $\gamma$ as a function of $\mu_{\lambda}$ and $C$.

Having shown that the performance index only depends on $\gamma$, we now define a probabilistic function $P\up{\gamma}_{\xi}(\gamma)$ of $\gamma$, which means that a battery with a normalized regulation energy capacity of $\gamma$ is $\xi$ certain to reach a performance score of  $P\up{\gamma}_{\xi}(\gamma)$. Hence we can rewrite the performance chance-constraint as 
\begin{align}
    \mathbf{Prob}\big[P(C\bm{r}, \bm{b}^g) \geq P\up{\gamma}_{\xi}(\gamma)\big] = \xi
\end{align}
and $P\up{\gamma}_{\xi}(\gamma)$ can be determined by simulating historical regulation signals assuming that the regulation signal distribution is stationary.
\begin{lemma}
$P\up{\gamma}_{\xi}(\gamma)$ is monotonic over $P\up{\gamma}_{\xi}(\gamma) \in (1-\delta, 1)$ and has an inverse function $[P\up{\gamma}_{\xi}]^{-1}(\rho)$ over $\rho\in (1-\delta, 1)$.
\end{lemma}
\emph{Proof:}
This lemma is trivial. First we consider a regulation signal realization set $\bm{r}$, the minimum performance index that a battery can possibly score is $1-\delta$ because the rest is not dependent on the battery response. Before reaching the perfect performance of $1$, an increment in $\gamma$ must result in an improvement in the performance index due to that $\gamma$ is the only constraining factor. This is true for any realizations of $\bm{r}$, hence it is trivial that $P\up{\gamma}_{\xi}(\gamma)$ is monotonic over $P\up{\gamma}_{\xi}(\gamma) \in (1-\delta, 1)$.

Following Lemma~1, if a participant wishes to reach a performance score $\rho\in (1-\delta, 1)$ with exactly $\xi$ confidence, it must use the following value of $\gamma$:
\begin{align}
    \gamma = [P\up{\gamma}_{\xi}]^{-1}(\rho),\quad \rho\in(1-\delta,1)
\end{align}
We now combine Theorem~2, \eqref{eq:pi_3} and \eqref{eq:pol3} to characterize analytically the optimal regulation capacity when using the proposed optimal response policy: 
\begin{align}\label{eq:optC1}
    C^*(\mu_{\lambda})
    &=\min\Bigg\{B, \frac{\min\{\overline{E}-\underline{E}, \hat{u}E\}}{[P\up{\gamma}_{\xi}]^{-1}(\rho\up{min})}\Bigg\}\nonumber\\
    &\text{where }
    \hat{u} =  \varphi^{-1}\Bigg(\frac{\eta^2+1}{\eta R \mu_{r}M}\delta\mu_{\lambda}\Bigg)
\end{align}
The only variable in this equation is $\mu_{\lambda}$. The rests are either based on the market policy or the BES design.
It is easy to see that if $\Phi$ is strictly convex, $C^*(\mu_{\lambda})$ is monotonic over the following range 
\begin{align}\label{eq:optC2}
    C&\in (0, \overline{C})\,\nonumber\\
    &\text{where } \overline{C} = \min\Big\{B, \frac{\overline{E}-\underline{E}}{[P\up{\gamma}_{\xi}]^{-1}(\rho\up{min})}\Big\}
\end{align}
and $\overline{C}$ is the maximum regulation capacity a participant could provide in order to reach the performance index $\rho\up{min}$ with a confidence of at least $\xi$. 

\subsection{Optimal Bidding Policy in Regulation Markets}

\begin{algorithm}[!htb]\label{alg:bid}
\SetAlgoLined
\KwResult{Determine the regulation offer price $\bm{\lambda}\up{b}$ associated with capacity segment $\bm{C}\up{b}$}
\tcp{\footnotesize starts from the first (cheapest) bid segment}
set $1 \to J$ 
\;
\tcp{\footnotesize goes through each segment until reaching $\overline{C}$}

\While{$\sum_{j=1}^J C\up{b}_j \leq \overline{C}$}{
  \tcp{\footnotesize total capacity offered so far}
  set $\sum_{j=1}^J C\up{b}_j \leq \overline{C} \to C\up{total}$\;
  \tcp{\footnotesize total payment price expected}
  set $[C^*]^{-1}(C\up{total}) \to \lambda\up{total}$\;
  \tcp{\footnotesize calculate the segment offer price}
  set $(\lambda\up{total}C\up{total}-\sum_{j=1}^{J-1} \lambda\up{b}_jC\up{b}_j)/C\up{b}_J \to \lambda\up{b}_J$\;
  \tcp{\footnotesize go to the next capacity segment}
  set $J+1 \to J$\;
 }
\caption{Optimal Battery Regulation Bidding Policy}
\end{algorithm}

At the bidding stage in the regulation market, each participant must submit a segment bidding curve that specifies how much regulation capacity that it is willing to provide at a given market price. We denote this set $(\bm{\lambda}\up{b}, \bm{C}\up{b})=\{(\lambda\up{b}_j, C\up{b}_j)\in \mathbb{R}_+^2\,|\,j\in \mathbb{N}\}$. Therefore, each participant needs a bidding policy to calculate $\bm{\lambda}\up{b}$ with respect to the regulation capacity segments $\bm{C}\up{b}$.

The optimal bidding policy is straightforwardly based on \eqref{eq:optC1} and \eqref{eq:optC2}: a participant uses the inverse function of $C^*(\cdot)$ to calculate the offer price associated with each capacity segment, while the total offered capacity must be smaller than $\overline{C}$ in order to satisfy the performance requirement. This optimal bidding policy is described in Algorithm~\ref{alg:bid}.


\section{Simulation}

\subsection{Data and Setting}

We use the following parameters for the battery energy storage in simulations unless otherwise specified:
\begin{itemize}
    \item Charging and discharging power rating: 10 MW
    \item Energy capacity: 3 MWh
    \item Charging and discharging efficiency: 95\%
    \item Maximum state of charge: 95\%
    \item Minimum state of charge: 10\%
    \item Round-trip efficiency: 92\%
    \item Battery cycle life: 1000 cycles at 80\% cycle depth
    \item Battery shelf life: 10 years
    \item Cell temperature: maintained at 25$^\circ C$
    \item Battery pack replacement cost: \$ 300~/kWh
    \item $\mathrm{Li(NiMnCo)O_2}$-based 18650 lithium-ion battery cells
\end{itemize}
These cells have a near-quadratic cycle depth stress function~\cite{laresgoiti2015modeling}:
\begin{align}
    \Phi(u) = (1.57\text{E-3})u\up{2.03}\,.
    \label{Eq:DoD}
\end{align}
We use a simple control policy as a benchmark in all of the following simulations. The simple policy uses all battery energy capacity $[\underline{E}, \overline{E}]$ for operation and does not consider the penalty price or the cycle aging cost, hence it is equivalent to fixing $\hat{u}=1$ in the proposed policy. Although this a naive control approach, it has been used in most stochastic real-time battery operation studies~\cite{oudalov2007optimizing,bitar2011role,zhang2017data}.

\subsection{Illustration of Theorem 1}\label{sec:sim_th1}

\begin{table}[!htb]
    \centering
    \caption{Simulation with randomly generated dispatch signals.}
    \begin{tabular}{l c c c c l l}
        \hline
        \hline
         Case & $\pi$ & $\eta$ & $T$ & $\hat{u}$  & \multicolumn{2}{c}{ Maximum regret [\$]} 
        \Tstrut\\
        & [\$/MWh] & [\%] & & [\%] & Proposed & Simple \Bstrut\\
        \hline
        1 & 50 & 100 & 100 & 11.1 &  {0.00}  & 183.9 \Tstrut\\
        2 & 100 & 100 & 100 & 21.9 &  {0.00} & 127.5 \\
        3 & 200 & 100 & 100 & 42.8 &  {0.00} & 47.9  \\
        4 & 50 & 92 & 100 & 11.2 &  {0.06}   & 184.9 \\
        5 & 50 & 92 & 200 & 11.2 &  {0.06}   & 408.8 \Bstrut\\
         \hline
    \end{tabular}
    \label{tab:sim}
\end{table}

We illustrate Theorem~1 by comparing the results of the proposed control policy with those of the simple control policy, Table~\ref{tab:sim} summarizes the simulation results. In each case the battery was subjected to 100 simulated regulation signals. Each case has a different simulation duration $T$, or penalty price $\pi$, or efficiency $\eta$. Control regrets are calculated by comparing the results from the control policy with the offline optimal control where all regulation signals are known. The proposed control policy achieves negligible regrets in all cases, while the simple policy control result has a significantly higher regrets. A further discussion of these results can be found in~\cite{xu2017optimal}.

\begin{table*}[!t]
    \centering
    \caption{Case Study Results in PJM RegD Market 03/2016-02/2017.}
    \begin{tabular}{r c c c c c c c}
        \hline
        \hline
                & Benchmark & \multicolumn{6}{c}{Optimal participation under $\xi$ performance confidence} \Tstrut\\
                &  & $\xi = 99\%$  & $\xi = 95\%$ & $\xi = 90\%$  & $\xi = 85\%$ & $\xi = 75\%$ & $\xi = 50\%$\Bstrut\\
        \hline
        Market income [k\$] & 1472.9 & 494.8 & 823.2 & 958.8 & 1035.9 & 1164.9 & 1266.8\Tstrut\Bstrut\\
        Aging cost [k\$] & 1091.8 & 67.5 & 165.8 & 213.0 & 240.3 & 286 & 326.3\Bstrut\\
        Prorated operating profit [k\$] &  381.1 & 427.3 & 657.4 & 745.7 & 795.6 & 878.6 & 940.5\Bstrut\\
        Battery cell life expectancy [month] & 9 & 69 & 42 & 35 & 33 & 29 & 26\Bstrut\\
        Annual average performance &  0.99 & 0.99 & 0.96 & 0.95 & 0.93 & 0.90 & 0.84
        \Bstrut\\
        Hours of under-performance & 6 & 81 & 204 & 246 & 297 & 426 & 1022
        \Bstrut\\
        Total regulation capacity cleared [MW$\cdot$h] & 87600 & 14991 & 29051 & 36838 & 42166 & 54087 & 73504
        \Bstrut\\
         \hline
    \end{tabular}
    \label{tab:case}
\end{table*}



\subsection{Illustration of Theorem 2}

\begin{figure}[!htb]%
	\centering
	\subfloat[]{
		\includegraphics[trim = 0mm 0mm 0mm 0mm, clip, width = .95\columnwidth]{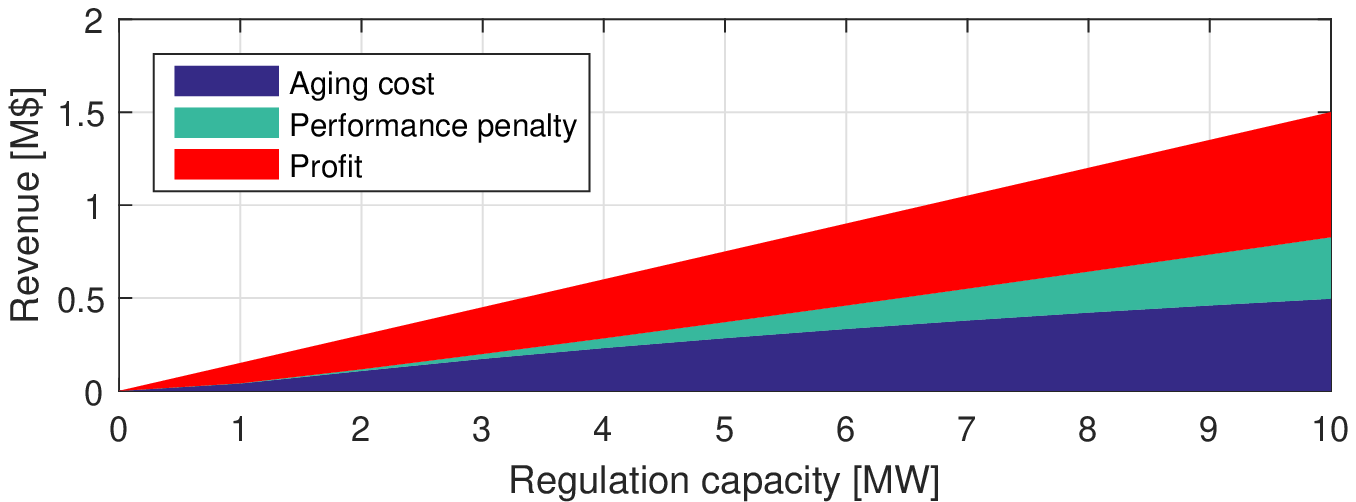}
		\label{fig:th2A}%
	}
	\\
	\subfloat[]{
		\includegraphics[trim = 0mm 0mm 0mm 0mm, clip, width = .95\columnwidth]{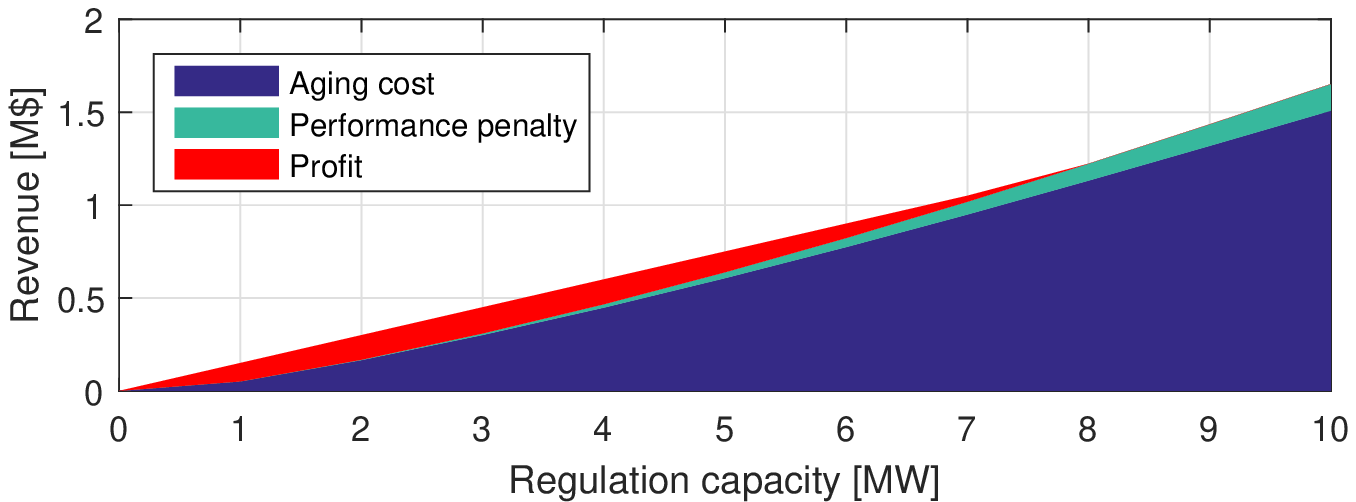}
		\label{fig:th2B}%
	}
	\caption{\footnotesize Full year revenue analysis of a 1~MWh BES participating in the regulation market with different regulation capacities. The total colored area reflects the payment a participant receives for regulation capacity. The operating profit is calculated by subtracting the aging and penalty costs from the total payment. (a) With the proposed control policy; (b) With the simple control policy.}%
	\label{Fig:th2}
\end{figure}

We illustrate Theorem~2 by simulating a 1~MWh BES providing a full-year regulation service with varying regulation capacities in the PJM market using price and signal data from 03/2016 to 02/2017. The minimum performance requirement is not enforced in this example. Fig.~\ref{Fig:th2} shows the simulation result for the proposed policy and the simple policy to demonstrate the effectiveness of our method. 

The result in Fig.~\ref{Fig:th2} includes the aging cost, the performance penalty cost, and the operational profit for regulation capacities up to 10~MW. With the proposed policy, the profits increase with the regulation capacity, which validates Theorem~2. With the simple policy, the BES operation is somewhat profitable when the regulation capacity is below 4~MW, but this profitability disappears completely once the regulation capacity is greater than 8~MW.

\subsection{Case Study}\label{sec:case}

We simulate a BES providing dynamic regulation (RegD) service in the PJM regulation market using signal and price data for a full-year (from 03/2016 to 02/2017) assuming the BES is a price taker and its bids has no effect on the market clearing price. The minimum performance index requirement for this market is 0.7, and the settlement period in one hour. $\delta$ in the performance index model \eqref{eq:lin_per} is set to 2/3 which approximates the PJM performance index calculation, this approximation is described in Appendix~C. All ex-post performance index calculations are performed according to the PJM manual~\cite{pjm_manual12}. 

\begin{figure}[!t]
\centering
	\subfloat[]{
		\includegraphics[trim = 0mm 0mm 0mm 0mm, clip, width = .95\columnwidth]{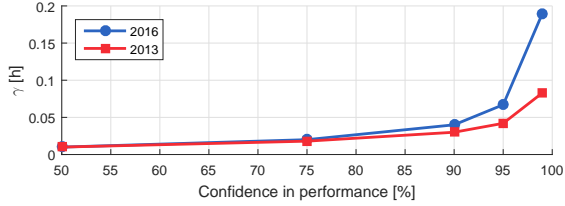}
		\label{fig:caseA}%
	}
	\\
	\subfloat[]{
		\includegraphics[trim = 0mm 0mm 0mm 0mm, clip, width = .95\columnwidth]{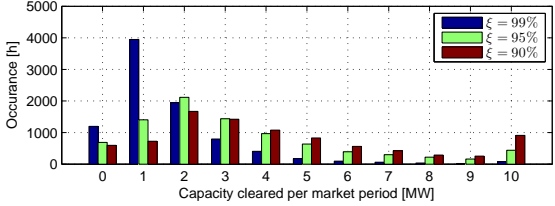}
		\label{fig:caseB}%
	}
    \caption{\footnotesize PJM case study. (a) The relative regulation capacity $\gamma$ vs. different $\xi$  from simulating historical regulation signals in 2013 and 2016; (b) Histogram of all cleared regulation capacities in 2016 from the proposed bidding policy under difference confidence $\xi$.}
\end{figure}

The regulation signal is slightly biased to make it energy zero-mean, a realistic assumption in the current market structure~\cite{xu2014bess}. Parameters of the bidding policy, including the value of $\mu_{r}$ and the function $P\up{\gamma}_{\xi}$, are determined using the PJM RegD signal from 06/2013 to 05/2014. The value of $\mu_{\lambda}$ is calculated from the PJM capacity and performance clearing prices using a mileage ratio of 3 based on the PJM market clearing manual~\cite{xu2016comparison}. We consider 10 1~MW regulation capacity offer segments for a 10MW/3MWh BES. The BES thus offers a total of 10~MW regulation capacity. The offer price associated with each segment is calculated according to the proposed bidding policy using the performance index simulation result from 2013, as shown in Fig.~\ref{fig:caseA}. Fig.~\ref{fig:caseB} shows the histogram of the regulation capacity cleared from the proposed bidding policy under different performance confidence. The performance confidence function concluded from the 2016 data is slightly different from 2013, in particular, the energy capacity requirement at high confidence levels are higher, possibly due to changes in the dispatch algorithm, unexpected weather conditions, or increased renewable penetration.

For comparison, we include a benchmark strategy, where the 10MW/3MWh BES participates the regulation market as a price taker and always offers 10~MW regulation capacity. Hence over the course of a year this BES provides a total of 87600~MW$\cdot$h of regulation capacity. The benchmark case uses a simple control strategy that always dispatches the entire BES energy capacity in response to regulation instructions.

\begin{figure}[!t]
\centering
	\subfloat[]{
		\includegraphics[trim = 0mm 0mm 0mm 0mm, clip, width = .95\columnwidth]{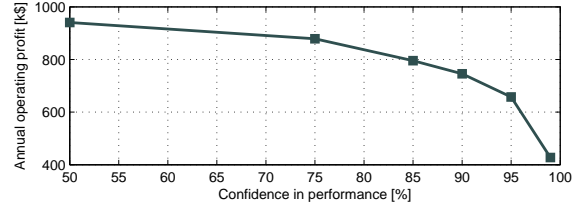}
		\label{fig:caseC}%
	}
	\\
	\subfloat[]{
		\includegraphics[trim = 0mm 0mm 0mm 0mm, clip, width = .95\columnwidth]{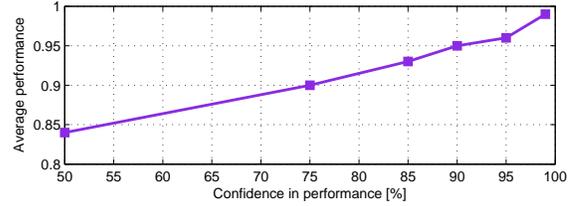}
		\label{fig:caseD}%
	}
    \caption{\footnotesize PJM case study risk analysis. (a) Annual prorated operating profit vs. performance confidence $\xi$; (b) Average annual performance index vs. performance confidence $\xi$.}
\end{figure}

Table~\ref{tab:case} summarizes the case study results, showing that the proposed policy earns more profit than the benchmark strategy for all performance confidence levels. On the other hand, the benchmark strategy has a higher regulation performance because it always dispatches the entire battery energy capacity. The operating profit increases with a lower performance confidence, i.e., a riskier strategy earns more profit. The downside of a riskier strategy is that the BES's offer may become noncompetitive because the system operator considers historical performance when clearing the market. We plot the operating profit (Fig.~\ref{fig:caseC}) and the average performance (Fig.~\ref{fig:caseD}) to show the trade-off between profit and risk. It is clear that while the average performance index decreases approximately linearly with $\xi$, the marginal increase in operating profit becomes smaller, i.e., the profit  from a riskier regulation market participation saturates. Therefore, a participant must carefully balance profit and risk in actual market scenarios.


\section{Conclusion}\label{sec:con}
In this paper, we proposed an optimal control policy and an optimal bidding policy for battery energy storage participating in performance-based regulation market. The proposed policies consider the cycle aging mechanism of electrochemical battery cells, and are adaptive to realistic market settings. We validate the optimality of the proposed policies using simulations, and demonstrate their effectiveness by a case study based on the PJM frequency regulation market. 

Performance-based frequency regulation markets have been proposed for around five years and their designs are still maturing. With an increasing number of battery participants, frequency regulation markets are becoming increasingly competitive. As exemplified by the ongoing market revisions in PJM, system operators may also tighten the performance requirements. We hope this study can inspire future researches on topics such as regulation market design considering profit-seeking battery participants, how will the regulation market saturates, and how will that impact the profitability of different types of participants.


\appendix

\subsection{Supplement Material to Theorem 1}
We calculate the optimal depth by transforming the optimization objective into a function of cycles instead of battery dispatch power $b_t$. We first observe that $b_t$ should never exceed the regulation instruction $Cr_t$ such that
\begin{align}
    0&\leq [b_t]^+ \leq C[r_t]^+ \label{eq:app01}\\
    0&\leq [-b_t]^+ \leq C[-r_t]^+ \label{eq:app02}
\end{align}
because otherwise the participant would undertake the performance penalty and the cost of degradation simultaneously. The mismatch term $||C\bm{r}-\bm{b}||_1$ in the objective function can be written as
\begin{align}\label{eq:app1}
    ||C\bm{r}-\bm{b}||_1 = C\mu_r - ||\bm{b}||_1
\end{align}
because $\mu_r = ||\bm{r}||_1$.
The next step is to represent $||\bm{b}||_1$ with a combination of cycles. We assume that the battery dispatch profile $\bm{b}$ consists a set of full cycles of depth $\bm{u}$ with no residue, i.e., all battery operation can be identified into cycles by the rainflow method. The summation all full cycles $\bm{u}$ must equal to the sum of all positive variations in the battery SoC as 
\begin{align}
    ||\bm{u}||_1 &= \frac{1}{E}\sum_{t=2}^T[e_t-e_{t-1}]^+ 
    =  \frac{M\eta}{E}\sum_{t=1}^T[b_t]^+
\end{align}
and similarly for negative SoC variations it follows that
\begin{align}
    ||\bm{u}||_1  &= \frac{1}{E}\sum_{t=2}^T[-e_t+e_{t-1}]^+ 
    =  \frac{M}{E\eta}\sum_{t=1}^T[-b_t]^+
\end{align}
and $||\bm{b}||_1$ can be rewritten as
\begin{align}\label{eq:app2}
    ||\bm{b}||_1 &= \sum_{t=1}^T[b_t]^+ + \sum_{t=1}^T[-b_t]^+
    = \frac{\eta^2+1}{\eta M}E||\bm{u}||_1\,.
\end{align}
We now substitute \eqref{eq:app1} and \eqref{eq:app2} into and rewrite the objective as a function of cycle depths (recall that all cycle depths are positive)
\begin{align}\label{eq:app03}
J(g,C) &= \mu_{\lambda} C - \mathbb{E}\Big[\delta\frac{\mu_{\lambda}}{\mu_{r}}
||C\bm{r}-\bm{b}^g||_1+A(\bm{b}^g)\Big] \\
&= (1-\delta)\mu_{\lambda} C 
 + \sum_{i=1}^{|\bm{u}|}\Big[\frac{\eta^2+1}{\eta}\pi E u_i-ER \Phi(u_i)\Big]\nonumber
\end{align}
where $\pi = \delta \mu_{\lambda}/(\mu_r M)$. \eqref{eq:app03} achieves the unconstrained minimum under the following condition
\begin{align}
    \frac{\partial J}{\partial u_i} &= 0, \qquad u_i = \hat{u} = \varphi^{-1}\Big(\frac{\eta^2+1}{\eta R}\pi\Big) \label{eq:app401}
\end{align}
where $\hat{u}$ is the optimal cycle depth, and the battery regulation operation is always profitable if all cycles are limited below $\hat{u}$. 
The proposed control policy therefore caps all cycle depths below $\hat{u}$ while enforcing operation constraints \eqref{eq:bat_con1}--\eqref{eq:bat_con3}, \eqref{eq:app01} and \eqref{eq:app02}. This optimal solution assumes all battery operation can be mapped into full cycles, however in most cases the rainflow method left a residue profile in which no more full cycles can be identified~\cite{rychlik1987new}. This leads to a control regret $\epsilon$ that is characterized by Theorem~1. $\epsilon$ is calculated as
\begin{align}
    \epsilon = 2J\ud{w}(\hat{u})+J\ud{v}(\hat{u})-2J\ud{w}(\hat{w})-J\ud{v}(\hat{v})
\end{align}
where
\begin{align}
    \hat{v} &= \varphi^{-1}\Big(\frac{\pi}{\eta R}\Big) \\
    \hat{w} &= \varphi^{-1}\Big(\frac{\pi\eta}{R}\Big) \\
    J\ud{v}(v) &= \frac{1}{\eta}\pi E v-\frac{1}{2}ER \Phi(v) \\
    J\ud{w}(w) &= \eta\pi E w -\frac{1}{2}ER \Phi(w) \,.
\end{align} 
Please refer to \cite{shi2017optimal} for further discussions and proof of this result.

\subsection{Proof for Theorem 2}

We start this proof by calculating the derivative of the response cycle depth with respect to the regulation capacity. In the proposed control policy, the battery follows regulation instruction until reaching $\hat{u}$, so the depth of a cycle scales linearly with $C$, then according to \eqref{eq:app2}
\begin{align}
    \frac{\partial u}{\partial C} = \begin{cases} 
    \eta M / [E (\eta^2+1)]  &  \text{if $u < \hat{u}$}\\
    0 & \text{else}\\
    \end{cases}
\end{align}
Then we calculate derivative of the expected profit with respect to $C$ based on \eqref{eq:app03}
\begin{align}
    &\frac{\partial J(g,C)}{\partial C} = (1-\delta)\mu_{\lambda} + \sum_{i=1}^{|\bm{u}|}\Big[\frac{\eta^2+1}{\eta}\pi E - ER \varphi(u_i)\Big]\frac{\partial u_i}{\partial C} \nonumber\\
    &= (1-\delta)\mu_{\lambda} + \frac{\eta M}{ \eta^2+1}\sum_{u_i < \hat{u}} \Big[\frac{\eta^2+1}{\eta}\pi  - R\varphi(u_i) \Big]
\end{align}
Which is strictly greater than zero because $\delta < 1$ and for $u_i < \hat{u}$
\begin{align}
    R \varphi(u_i) < R \varphi(\hat{u}) = (\eta^2+1)\pi/\eta\,,
\end{align}
when $\Phi(\cdot)$ is strictly convex and $\varphi = d\Phi/d u$ is an strictly increasing function. Therefore, the optimal regulation capacity $C^{*}$ that maximizes the expected operating profit is equal to maximum possible capacity bidding, which is the BES power rating $B$. Theorem~2 is hence proved.

\subsection{PJM Performance Index Simplification}
PJM penalizes the entire regulation payment with the performance index, which is calculated with equal weights on precision, correlation, and delay~\cite{pjm_manual12}:
\begin{align}
    &\mathrm{Performance Index} =\\
    &\frac{1}{3}\Big(\mathrm{Precision Score} + \mathrm{Correlation Score} + \mathrm{Delay Score}\Big)\nonumber
\end{align}
while all scores are valued between 0 to 1. The precision score is the absolute error of the regulation response as
\begin{align}
    \mathrm{Precision} = 1-\frac{||C\bm{r}-\bm{b}||_1}{||C\bm{r}||_1}
\end{align}
The delay score is evaluated over the delay under which the time-shifted BES response achieves the highest correlation with the regulation instruction, and this correlation coefficient is the correlation score. Because BES responses to dispatch commands instantaneously, the response profile always achieves highest correlation with the instruction at zero delay, hence the BES always has full delay score. However considering the difficulty of incorporating correlation index calculation into a convex programming problem, we approximate the correlation score to be the same as the precision score. The PJM performance index calculation is therefore simplified to
\begin{align}
    \mathrm{Performance Index} &=\frac{1}{3}\Big(2-2\frac{||C\bm{r}-\bm{b}||_1}{||C\bm{r}||_1}+1\Big) \nonumber\\
    &= 1 - \frac{2}{3}\frac{||C\bm{r}-\bm{b}||_1}{||C\bm{r}||_1}
\end{align}
which is the same as setting $\delta=2/3$ in \eqref{eq:lin_per}.

\bibliographystyle{IEEEtran}	
\bibliography{IEEEabrv,literature}		

\begin{thebibliography}{10}
\providecommand{\url}[1]{#1}
\csname url@samestyle\endcsname
\providecommand{\newblock}{\relax}
\providecommand{\bibinfo}[2]{#2}
\providecommand{\BIBentrySTDinterwordspacing}{\spaceskip=0pt\relax}
\providecommand{\BIBentryALTinterwordstretchfactor}{4}
\providecommand{\BIBentryALTinterwordspacing}{\spaceskip=\fontdimen2\font plus
\BIBentryALTinterwordstretchfactor\fontdimen3\font minus
  \fontdimen4\font\relax}
\providecommand{\BIBforeignlanguage}[2]{{%
\expandafter\ifx\csname l@#1\endcsname\relax
\typeout{** WARNING: IEEEtran.bst: No hyphenation pattern has been}%
\typeout{** loaded for the language `#1'. Using the pattern for}%
\typeout{** the default language instead.}%
\else
\language=\csname l@#1\endcsname
\fi
#2}}
\providecommand{\BIBdecl}{\relax}
\BIBdecl

\bibitem{irena}
\BIBentryALTinterwordspacing
{The International Renewable Energy Agency (IRENA) }, ``Battery storage for
  renewables: Market status and technology outlook.'' [Online]. Available:
  \url{http://www.irena.org/documentdownloads/publications/irena_battery_storage_report_2015.pdf}
\BIBentrySTDinterwordspacing

\bibitem{pjm2016}
\BIBentryALTinterwordspacing
``2016 state of the market report for pjm.'' [Online]. Available:
  \url{http://www.monitoringanalytics.com/reports/PJM_State_of_the_Market/2016.shtml}
\BIBentrySTDinterwordspacing

\bibitem{nykvist2015rapidly}
B.~Nykvist and M.~Nilsson, ``Rapidly falling costs of battery packs for
  electric vehicles,'' \emph{Nature Climate Change}, vol.~5, no.~4, pp.
  329--332, 2015.

\bibitem{qiu2016stochastic}
T.~Qiu, B.~Xu, Y.~Wang, Y.~Dvorkin, and D.~Kirschen, ``Stochastic multi-stage
  co-planning of transmission expansion and energy storage,'' \emph{IEEE
  Transactions on Power Systems}, vol.~PP, no.~99, pp. 1--1, 2016.

\bibitem{FERC755}
\BIBentryALTinterwordspacing
FERC, ``Frequency regulation compensation in the organized wholesale power
  markets (order no. 755).'' [Online]. Available:
  \url{https://www.ferc.gov/whats-new/comm-meet/2011/102011/E-28.pdf}
\BIBentrySTDinterwordspacing

\bibitem{FERC784}
\BIBentryALTinterwordspacing
``Third-party provision of ancillary services; accounting and financial
  reporting for new electric storage technologies (order no. 784).'' [Online].
  Available:
  \url{https://www.ferc.gov/whats-new/comm-meet/2013/071813/E-22.pdf}
\BIBentrySTDinterwordspacing

\bibitem{xu2016comparison}
B.~Xu, Y.~Dvorkin, D.~S. Kirschen, C.~Silva-Monroy, and J.-P. Watson, ``A
  comparison of policies on the participation of storage in us frequency
  regulation markets,'' \emph{PES Genreal Meeting}, 2016.

\bibitem{xu2014bess}
B.~Xu, A.~Oudalov, J.~Poland, A.~Ulbig, and G.~Andersson, ``Bess control
  strategies for participating in grid frequency regulation,'' in \emph{19th
  IFAC World Congress, Cape Town}, 2014.

\bibitem{mercier2009optimizing}
P.~Mercier, R.~Cherkaoui, and A.~Oudalov, ``Optimizing a battery energy storage
  system for frequency control application in an isolated power system,''
  \emph{IEEE Transactions on Power Systems}, vol.~24, no.~3, pp. 1469--1477,
  2009.

\bibitem{borsche2013power}
T.~Borsche, A.~Ulbig, M.~Koller, and G.~Andersson, ``Power and energy capacity
  requirements of storages providing frequency control reserves,'' in
  \emph{IEEE PES General Meeting, Vancouver}, 2013.

\bibitem{watson2017comparing}
D.~Watson, C.~Hastie, and M.~Rodgers, ``Comparing different regulation
  offerings from a battery in a wind r d park,'' \emph{IEEE Transactions on
  Power Systems}, vol.~PP, no.~99, pp. 1--1, 2017.

\bibitem{tu2017maximizing}
T.~A. Nguyen, R.~H. Byrne, R.~J. Concepcion, and I.~Gyuk, ``Maximizing revenue
  from electrical energy storage in miso energy \& frequency regulation
  markets,'' in \emph{IEEE PES General Meeting, Chicago}, 2017.

\bibitem{byrne2012estimating}
R.~H. Byrne and C.~A. Silva-Monroy, ``Estimating the maximum potential revenue
  for grid connected electricity storage: Arbitrage and regulation,''
  \emph{Sandia National Laboratories}, 2012.

\bibitem{zhang2017data}
H.~Zhang, Z.~Hu, E.~Munsing, S.~J. Moura, and Y.~Song, ``Data-driven
  chance-constrained regulation capacity offering for distributed energy
  resources,'' \emph{arXiv preprint arXiv:1708.05114}, 2017.

\bibitem{swierczynski2014selection}
M.~Swierczynski, D.~I. Stroe, A.-I. Stan, R.~Teodorescu, and D.~U. Sauer,
  ``{Selection and performance-degradation modeling of
  LiMO$_2$/Li$_4$Ti$_5$O$_{12}$ and LiFePO$_4$/C battery cells as suitable
  energy storage systems for grid integration with wind power plants: an
  example for the primary frequency regulation service},'' \emph{IEEE
  transactions on Sustainable Energy}, vol.~5, no.~1, pp. 90--101, 2014.

\bibitem{stroe2015degradation}
D.-I. Stroe, M.~Swierczynski, A.-I. Stroe, R.~Teodorescu, R.~Laerke, and P.~C.
  Kjaer, ``Degradation behaviour of lithium-ion batteries based on field
  measured frequency regulation mission profile,'' in \emph{Energy Conversion
  Congress and Exposition (ECCE), 2015 IEEE}.\hskip 1em plus 0.5em minus
  0.4em\relax IEEE, 2015, pp. 14--21.

\bibitem{koller2015review}
M.~Koller, T.~Borsche, A.~Ulbig, and G.~Andersson, ``Review of grid
  applications with the zurich 1mw battery energy storage system,''
  \emph{Electric Power Systems Research}, vol. 120, pp. 128--135, 2015.

\bibitem{white2011using}
C.~D. White and K.~M. Zhang, ``Using vehicle-to-grid technology for frequency
  regulation and peak-load reduction,'' \emph{Journal of Power Sources}, vol.
  196, no.~8, pp. 3972--3980, 2011.

\bibitem{stroe2016degradation}
D.-I. Stroe, M.~Swierczynski, A.-I. Stroe, R.~Laerke, P.~C. Kjaer, and
  R.~Teodorescu, ``Degradation behavior of lithium-ion batteries based on
  lifetime models and field measured frequency regulation mission profile,''
  \emph{IEEE Transactions on Industry Applications}, vol.~52, no.~6, pp.
  5009--5018, 2016.

\bibitem{zhang2016profit}
Y.~J. Zhang, C.~Zhao, W.~Tang, and S.~H. Low, ``Profit maximizing planning and
  control of battery energy storage systems for primary frequency control,''
  \emph{IEEE Transactions on Smart Grid}, 2016.

\bibitem{cui2015optimal}
T.~Cui, Y.~Wang, S.~Chen, Q.~Zhu, S.~Nazarian, and M.~Pedram, ``Optimal control
  of pevs for energy cost minimization and frequency regulation in the smart
  grid accounting for battery state-of-health degradation,'' in
  \emph{Proceedings of the 52nd Annual Design Automation Conference}.\hskip 1em
  plus 0.5em minus 0.4em\relax ACM, 2015, p. 134.

\bibitem{he2015optimal}
G.~He, Q.~Chen, C.~Kang, P.~Pinson, and Q.~Xia, ``Optimal bidding strategy of
  battery storage in power markets considering performance-based regulation and
  battery cycle life,'' \emph{IEEE Transactions on Smart Grid}, vol.~7, no.~5,
  pp. 2359--2367, Sept 2016.

\bibitem{pjm_manual12}
\BIBentryALTinterwordspacing
``Pjm manual 12: Balancing operations.'' [Online]. Available:
  \url{http://www.pjm.com/~/media/documents/manuals/m12.ashx}
\BIBentrySTDinterwordspacing

\bibitem{ramadesigan2012modeling}
V.~Ramadesigan, P.~W. Northrop, S.~De, S.~Santhanagopalan, R.~D. Braatz, and
  V.~R. Subramanian, ``Modeling and simulation of lithium-ion batteries from a
  systems engineering perspective,'' \emph{Journal of The Electrochemical
  Society}, vol. 159, no.~3, pp. R31--R45, 2012.

\bibitem{ning2004cycle}
G.~Ning and B.~N. Popov, ``Cycle life modeling of lithium-ion batteries,''
  \emph{Journal of The Electrochemical Society}, vol. 151, no.~10, pp.
  A1584--A1591, 2004.

\bibitem{reniers2018improving}
\BIBentryALTinterwordspacing
J.~M. Reniers, G.~Mulder, S.~Ober-Bl{\"o}baum, and D.~A. Howey, ``Improving
  optimal control of grid-connected lithium-ion batteries through more accurate
  battery and degradation modelling,'' \emph{Journal of Power Sources}, vol.
  379, pp. 91 -- 102, 2018. [Online]. Available:
  \url{http://www.sciencedirect.com/science/article/pii/S0378775318300041}
\BIBentrySTDinterwordspacing

\bibitem{ying2016stochastic}
W.~Ying, Z.~Zhi, A.~Botterud, K.~Zhang, and D.~Qia, ``Stochastic coordinated
  operation of wind and battery energy storage system considering battery
  degradation,'' \emph{Journal of Modern Power Systems and Clean Energy},
  vol.~4, no.~4, pp. 581--592, 2016.

\bibitem{wang2017improving}
Y.~Wang, C.~Wan, Z.~Zhou, K.~Zhang, and A.~Botterud, ``Improving deployment
  availability of energy storage with data-driven agc signal models,''
  \emph{IEEE Transactions on Power Systems}, 2017.

\bibitem{megel2013maximizing}
O.~Megel, J.~L. Mathieu, and G.~Andersson, ``Maximizing the potential of energy
  storage to provide fast frequency control,'' in \emph{Innovative Smart Grid
  Technologies Europe (ISGT EUROPE), 2013 4th IEEE/PES}.\hskip 1em plus 0.5em
  minus 0.4em\relax IEEE, 2013, pp. 1--5.

\bibitem{fortenbacher2017modeling}
P.~Fortenbacher, J.~L. Mathieu, and G.~Andersson, ``Modeling and optimal
  operation of distributed battery storage in low voltage grids,'' \emph{IEEE
  Transactions on Power Systems}, vol.~PP, no.~99, pp. 1--1, 2017.

\bibitem{ortega2014optimal}
M.~A. Ortega-Vazquez, ``Optimal scheduling of electric vehicle charging and
  vehicle-to-grid services at household level including battery degradation and
  price uncertainty,'' \emph{IET Generation, Transmission \& Distribution},
  vol.~8, no.~6, pp. 1007--1016, 2014.

\bibitem{koller2013defining}
M.~Koller, T.~Borsche, A.~Ulbig, and G.~Andersson, ``Defining a degradation
  cost function for optimal control of a battery energy storage system,'' in
  \emph{PowerTech (POWERTECH), 2013 IEEE Grenoble}.\hskip 1em plus 0.5em minus
  0.4em\relax IEEE, 2013, pp. 1--6.

\bibitem{wang2016quantifying}
D.~Wang, J.~Coignard, T.~Zeng, C.~Zhang, and S.~Saxena, ``Quantifying electric
  vehicle battery degradation from driving vs. vehicle-to-grid services,''
  \emph{Journal of Power Sources}, vol. 332, pp. 193--203, 2016.

\bibitem{shi2017using}
Y.~Shi, B.~Xu, D.~Wang, and B.~Zhang, ``Using battery storage for peak shaving
  and frequency regulation: Joint optimization for superlinear gains,''
  \emph{IEEE Transactions on Power Systems}, 2017.

\bibitem{nguyen2016stochastic}
T.~A. Nguyen and M.~Crow, ``Stochastic optimization of renewable-based
  microgrid operation incorporating battery operating cost,'' \emph{IEEE
  Transactions on Power Systems}, vol.~31, no.~3, pp. 2289--2296, 2016.

\bibitem{ecker2014calendar}
M.~Ecker, N.~Nieto, S.~K{\"a}bitz, J.~Schmalstieg, H.~Blanke, A.~Warnecke, and
  D.~U. Sauer, ``Calendar and cycle life study of li (nimnco) o 2-based 18650
  lithium-ion batteries,'' \emph{Journal of Power Sources}, vol. 248, pp.
  839--851, 2014.

\bibitem{ruetschi2004aging}
P.~Ruetschi, ``Aging mechanisms and service life of lead--acid batteries,''
  \emph{Journal of Power Sources}, vol. 127, no.~1, pp. 33--44, 2004.

\bibitem{xu2016modeling}
B.~Xu, A.~Oudalov, A.~Ulbig, G.~Andersson, and D.~Kirschen, ``Modeling of
  lithium-ion battery degradation for cell life assessment,'' \emph{IEEE
  Transactions on Smart Grid}, vol.~PP, no.~99, pp. 1--1, 2016.

\bibitem{wang2014degradation}
J.~Wang, J.~Purewal, P.~Liu, J.~Hicks-Garner, S.~Soukazian, E.~Sherman,
  A.~Sorenson, L.~Vu, H.~Tataria, and M.~W. Verbrugge, ``Degradation of lithium
  ion batteries employing graphite negatives and nickel--cobalt--manganese
  oxide+ spinel manganese oxide positives: Part 1, aging mechanisms and life
  estimation,'' \emph{Journal of Power Sources}, vol. 269, pp. 937--948, 2014.

\bibitem{rychlik1987new}
I.~Rychlik, ``A new definition of the rainflow cycle counting method,''
  \emph{International journal of fatigue}, vol.~9, no.~2, pp. 119--121, 1987.

\bibitem{abdulla2016optimal}
K.~Abdulla, J.~De~Hoog, V.~Muenzel, F.~Suits, K.~Steer, A.~Wirth, and
  S.~Halgamuge, ``Optimal operation of energy storage systems considering
  forecasts and battery degradation,'' \emph{IEEE Transactions on Smart Grid},
  2016.

\bibitem{jin2018applicability}
X.~Jin, A.~Vora, V.~Hoshing, T.~Saha, G.~Shaver, O.~Wasynczuk, and
  S.~Varigonda, ``Applicability of available li-ion battery degradation models
  for system and control algorithm design,'' \emph{Control Engineering
  Practice}, vol.~71, pp. 1--9, 2018.

\bibitem{laresgoiti2015modeling}
I.~Laresgoiti, S.~K{\"a}bitz, M.~Ecker, and D.~U. Sauer, ``Modeling mechanical
  degradation in lithium ion batteries during cycling: Solid electrolyte
  interphase fracture,'' \emph{Journal of Power Sources}, vol. 300, pp.
  112--122, 2015.

\bibitem{millner2010modeling}
A.~Millner, ``Modeling lithium ion battery degradation in electric vehicles,''
  in \emph{Innovative Technologies for an Efficient and Reliable Electricity
  Supply (CITRES), 2010 IEEE Conference on}.\hskip 1em plus 0.5em minus
  0.4em\relax IEEE, 2010, pp. 349--356.

\bibitem{omar2014lithium}
N.~Omar, M.~A. Monem, Y.~Firouz, J.~Salminen, J.~Smekens, O.~Hegazy,
  H.~Gaulous, G.~Mulder, P.~Van~den Bossche, T.~Coosemans \emph{et~al.},
  ``Lithium iron phosphate based battery--assessment of the aging parameters
  and development of cycle life model,'' \emph{Applied Energy}, vol. 113, pp.
  1575--1585, 2014.

\bibitem{bell1982regret}
D.~E. Bell, ``Regret in decision making under uncertainty,'' \emph{Operations
  research}, vol.~30, no.~5, pp. 961--981, 1982.

\bibitem{shi2017optimal}
Y.~Shi, B.~Xu, and B.~Zhang, ``Optimal battery control under cycle aging
  mechanisms,'' \emph{arXiv preprint arXiv:1709.05715}, 2017.

\bibitem{oudalov2007optimizing}
A.~Oudalov, D.~Chartouni, and C.~Ohler, ``Optimizing a battery energy storage
  system for primary frequency control,'' \emph{IEEE Transactions on Power
  Systems}, vol.~22, no.~3, pp. 1259--1266, 2007.

\bibitem{bitar2011role}
E.~Bitar, R.~Rajagopal, P.~Khargonekar, and K.~Poolla, ``The role of co-located
  storage for wind power producers in conventional electricity markets,'' in
  \emph{American Control Conference (ACC), 2011}.\hskip 1em plus 0.5em minus
  0.4em\relax IEEE, 2011, pp. 3886--3891.

\bibitem{xu2017optimal}
B.~Xu, Y.~Shi, D.~S. Kirschen, and B.~Zhang, ``Optimal regulation response of
  batteries under cycle aging mechanisms,'' in \emph{IEEE Conference on
  Decision and Control, Melbourne}, 2017.

\end{thebibliography}

\begin{IEEEbiographynophoto}{Bolun Xu}
received B.S. degrees in Electrical and Computer Engineering
from Shanghai Jiaotong
University, Shanghai, China in 2011, and the M.Sc degree in Electrical
Engineering from Swiss Federal Institute of Technology, Zurich, Switzerland
in 2014.

He is currently pursuing the Ph.D. degree in Electrical Engineering at the
University of Washington, Seattle, WA, USA. His research interests include
energy storage, power system operations, and power system economics.
\end{IEEEbiographynophoto}

\begin{IEEEbiographynophoto}{Yuanyuan Shi} 
received her Bachelor of Engineering in Automation from Nanjing Univeristy, China in 2015. She is currently pursuing her Ph.D. degree in Electrical Engineering and M.S. in Statistics at the University of Washington, Seattle, WA, USA.

She works on optimization, control and machine learning, with applications in power systems and cyberphysical systems. She is a recipient of Washington Clean Energy Institute graduate fellowship on energy and sustainability.
\end{IEEEbiographynophoto}

\begin{IEEEbiographynophoto}{Daniel S. Kirschen}
received his electrical and mechanical engineering degree from the Universite Libre de Bruxelles, Brussels, Belgium, in 1979 and his M.S. and Ph.D. degrees from the University of Wisconsin, Madison, WI, USA, in 1980, and 1985, respectively.

He is currently the Donald W. and Ruth Mary Close Professor of Electrical Engineering at the University of Washington, Seattle, WA, USA. His research interests include smart grids, the integration of renewable energy sources in the grid, power system economics, and power system security.
\end{IEEEbiographynophoto}

\begin{IEEEbiographynophoto}{Baosen Zhang}
(M'86-SM'91-F'07) received his Bachelor of Applied Science in Engineering Science degree from the University of Toronto in 2008; and his PhD degree in Electrical Engineering and Computer Sciences from University of California, Berkeley in 2013.

He was a Postdoctoral Scholar at Stanford University, affiliated with the Civil and Environmental Engineering and Management \& Science Engineering. He is currently an Assistant Professor in Electrcial Engineering at the University of Washington, Seattle, WA. His research interests are in power systems and cyberphysical systems. He was selected as one of Forbe's 30 under 30 in energy in 2015.
\end{IEEEbiographynophoto}

\end{document}